\documentclass[a4paper]{amsart}

\usepackage[dvips]{graphicx}

\usepackage{epsf,latexsym}

\begin{document}

\title{Automatic structures and growth functions \\
for finitely generated abelian groups}
\author{KAMEI Satoshi}
\date{}
\address{School of Computer Science, Tokyo University of Technology, 1404-1 Katakura, Hachioji-shi Tokyo 192-0982, Japan}

\email{kamei@cs.teu.ac.jp}


\maketitle

\begin{abstract}
In this paper, we consider the formal power series whose $n$-th coefficient is the number of copies of a given finite graph in the ball of radius $n$ centred at the identity element in the Cayley graph of a finitely generated group and call it the growth function. Epstein, Iano-Fletcher and Uri Zwick proved that the growth function is a rational function if the group has a geodesic automatic structure. 
We compute the growth function in the case where the group is abelian and see that the denominator of the rational function is determined from the rank of the group.

\end{abstract}

\section{Introduction}

The theory of automatic groups was introduced in \cite{CE}. If a group has an automatic structure, there exist efficient algorithms for enumerating elements and checking two elements which differ by multiplication on the right by only one generator. 
In~\cite{CE} and further works, it was shown that many classes of finitely generated groups, such as hyperbolic groups, finitely generated abelian groups, braid groups, have automatic structures. 
It is known that part of the definition of automatic structure can be replaced with a property of paths in Cayley graphs. Thus automatic structures for groups are closely related to geometric structures of their Cayley graphs.

In this paper, we observe a discrete analogue of the growth of volumes which are known in the field of Riemannian geometry. 
In short, for an automatic group, we consider the numbers of finite graphs isomorphic to a given finite graph in the balls of various radii centred at the identity element in the Cayley graph and compute the growth series of the numbers by use of a finite state automaton which is part of an automatic structure. 
First, we prepare some notations to describe our concern and the main result. 

Let $G$ be a finitely generated group. A generating set $\Sigma$ for $G$ is called {\it symmetric} if $x \in \Sigma$, then $x^{-1}$ is also in $\Sigma$. In the followings, we always assume that generating sets for groups are symmetric. 
Let $\Gamma(G, \Sigma)$ denote the Cayley graph of $(G, \Sigma)$. For simplicity, we usually write $\Gamma$ instead of $\Gamma(G, \Sigma)$. We consider that the length of each edge of $\Gamma$ is equal to $1$ and extend the distance to all points of $\Gamma$. This makes $\Gamma$ into a geodesic metric space. Let $\Gamma_n$ be the ball of radius $n$ in $\Gamma$ centred at the identity element of $G$. 

Let $S$ be a directed subgraph of $\Gamma$ whose edges are labelled by elements of $\Sigma$. A {\it morphism} $f:S\to \Gamma$ is defined to be a function that maps each vertex to a vertex and each labelled directed edge to a directed edge with the same label. Furthermore the initial and final endpoints of a directed edge in $S$ are mapped to the initial and final endpoints respectively of the corresponding directed edge of $\Gamma$. 
Let $b_n(S)$ be the number of morphisms from $S$ to $\Gamma_n$. Set $c_n(S) = b_n(S) -b_{n-1}(S)$ for $n \ge 1$ and $c_0(S) = b_0(S)$. 
Correspondingly we consider the formal power series $C(S, z) = \sum_{j =0}^{\infty} c_j(S) z^j$ and call it the {\it growth function} for $S$. 
In the case where $S$ is a vertex, $c_n(S)$ counts the number of elements of $G$ whose shortest representatives with respect to $(G, \Sigma)$ have exactly length $n$. 
In~\cite{Ch}, Charny computed the growth functions of Artin groups of finite type in that case.

The problem of computing the growth functions for various subgraphs arose from attempts to generalise the Ising model in quantum mechanics to more general geometries by Saito(\cite{S}). 
In~\cite{E}, Epstein, Iano-Fletcher and Uri Zwick computed growth functions in many cases. They proved that $C(S,z)$ is a rational function if $G$ has a geodesic automatic structure. Furthermore they showed that one can choose a common denominator of the rational function independently of $S$ if $G$ is hyperbolic and $\Sigma$ is fixed. 
They pointed out that similar assertions seem to hold in other classes of automatic groups. 
The aim of this paper is to consider the problem in the case where the groups are abelian. The main theorem is the following.

\vspace*{\baselineskip}

{\bf The main theorem.} {\it Let G be a finitely generated abelian group whose rank is $r$ and $\Sigma$ be a symmetric generating set. Let $S$ be any non-empty finite connected labelled directed subgraph of $\Gamma(G, \Sigma)$. Then $C(S, z)$ is a rational function with integral coefficients and the denominator is $(1-z)^r$.}

\vspace*{\baselineskip}

In section $2$, we recall the definition of automatic structure. Also we see that finitely generated abelian groups are automatic. In section $3$, we construct finite state automata which are part of automatic structures for finitely generated abelian groups. Furthermore we prepare some lemmas to compute the growth functions. In section $4$, we prove the main theorem. In section $5$, we see some simple examples.

\section{Automatic structure}

In the followings, a finitely generating set for a group is also called an {\it alphabet}. Each element of an alphabet is called a {\it letter}. A {\it word} over an alphabet is a finite sequence of letters of the alphabet. For an alphabet $\Sigma$, we denote by $\Sigma^{*}$ the set of all words over $\Sigma$. A subset of $\Sigma^*$ is called a {\it language} over $\Sigma$. A {\it regular language} is the set of all words accepted by a finite state automaton. See \cite{CE} for the basic notations about finite state automata. For a finite state automaton $A$, here and subsequently $L(A)$ denotes the regular language consists of all words accepted by $A$. There is a natural map $\pi : \Sigma^{*} \to G$ since $\Sigma^{*}$ can be seen as the free monoid by $\Sigma$. For a word $w \in \Sigma^{*}$, we denote its image in $G$ by $\overline{w}$ instead of $\pi (w)$ for simplicity.

\vspace*{\baselineskip}

{\sc Definition} (automatic structure).  Let $G$ be a group and $\Sigma$ be a symmetric generating set. $G$ is said to be {\it automatic} if there exist a finite state automaton $W$ over $\Sigma$ and a finite state automaton $M_x$ over $(\Sigma, \Sigma)$ for each $x \in \Sigma \cup \{ e \}$ which satisfy the following conditions.

\noindent
(1) The map $\pi :L(W) \to G$ is surjective.

\noindent
(2) The pair of strings $(w_1, w_2) \in (\Sigma, \Sigma)$ is accepted by $M_x$ if and only if $w_1$ and 

 \ $w_2$ are accepted by $W$ and $\overline{w_1x}= \overline{w_2}$.

\vspace*{\baselineskip}

We see a relation between automatic structures for groups and geometric structures of their Cayley graphs. For the purpose, we define a property of paths in geodesic metric spaces.

\vspace*{\baselineskip}
{\sc Definition} (fellow travellers property). Let $(X,d)$ be a geodesic metric space. Let $\alpha : [0, a] \to X$ and $\beta : [0,b] \to X$ be two paths (not necessarily geodesics). If $a \ge b$, we set $\beta(t) = \beta (b)$ for all $b \le t \le a$. Fix a constant $\kappa \ge 0$. We say $\alpha$ and $\beta$ are $\kappa$-{\it fellow travellers} if $d(\alpha(t), \beta(t))  \le \kappa$ for all $t$ such that $0 \le t \le \max (a,b)$. The constant $\kappa$ is called a {\it fellow travellers constant}.

\vspace*{\baselineskip}

The next theorem guarantees that fellow travellers property is equivalent to the condition (2) in the definition of automatic structure. For a word $w \in L(W)$, $\hat{w}$ denotes the path from the identity element to $\overline{w}$ in $\Gamma$ according to $w$.

\vspace*{\baselineskip}

{\bf Theorem 2.1 (\cite{CE}).} {\it Let $W$ be a finite state automaton and $L(W)$ be a regular language such that the map $\pi : L(W) \to G$ is surjective. Then $W$ is part of an automatic structure for $(G, \Sigma)$ if and only if there exists $\kappa \ge 1$ such that for $w_1, w_2 \in L(W)$ and $x \in \Sigma \cup \{e\} $ satisfying $\overline{w_1 x} =\overline{w_2}$, $\hat{w}_1$ and $\hat{w}_2$ are $\kappa$-fellow travellers in $\Gamma(G, \Sigma)$}. 

\vspace*{\baselineskip}

In the followings, we say that $L(W)$ satisfies $\kappa$-{\it fellow travellers property} if there exists a constant $\kappa$ stated in Theorem 2.1.

In the remain of this section, we see that finitely generated abelian groups have automatic structures. First, we define an order of words.

\vspace*{\baselineskip}
{\sc Definition} (shortlex order). Let $\Sigma$ be an ordered alphabet. 
For any pair of words $v, w \in \Sigma^{*}$, we define $v <w$ if $v$ is shorter than $w$, or they have the same length and $v$ comes before $w$ in lexicographical order. The order obtained in this manner is called the {\it shortlex order} of $\Sigma^{*}$.

\vspace*{\baselineskip}

A word $w \in \Sigma^{*}$ is called {\it shortlex} if $w$ is the minimum representative in the shortlex order among all words which represent the same element $\overline{w}$. If there exists a word acceptor $W$ for $(G, \Sigma)$ which only accepts shortlex words, we say that $G$ has a {\it shortlex automatic structure}.

Notice that if a word acceptor $W$ is part of a shortlex automatic structure, then the natural map $\pi : L(W) \rightarrow G$ is bijective and all paths in $\Gamma$ according to the words of L(W) are geodesics in $\Gamma(G, \Sigma)$. Further $L(W)$ is {\it prefix closed}, which means that all prefixes of any word in $L(W)$ are also accepted by $W$.

\vspace*{\baselineskip}
{\bf Theorem 2.2 (\cite[Theorem 4.3.1]{CE}).} {\it Let $G$ be a finitely generated abelian group and $\Sigma$ be a symmetric generating set. Then $G$ has a shortlex automatic structure with respect to any ordered set of elements of $\Sigma$}.

\vspace*{\baselineskip}

{\rm Sketch of Proof.} It is easily seen that there exists a finite state automaton which only accepts the shortlex words for $(G, \Sigma)$ with any order in $\Sigma$. 
Thus we only verify that the shortlex words satisfy the fellow travellers property. 

Let $\{a_1, \cdots, a_n\}$ be an ordered alphabet. We put every word in {\it normal form} $a_1^{r_1} \cdots a_n^{r_n}$ such that $r_1, \cdots , r_n$ are non-negative integers. Further we identify $a_1^{r_1} \cdots a_n^{r_n}$ with the $n$-tuple $(r_1, \cdots , r_n)$. For $v = a_1^{r_1} \cdots a_n^{r_n}$ and $w = a_1^{s_1} \cdots a_n^{s_n}$, we say that $v$ is contained in $w$ if $r_i \le s_i$ for all $i$ and write $v \prec w$.

If $w$ and $w'$ are distinct words which represent the same group element, we say that we have a {\it relation} $w \sim w'$. A relation $v \sim v'$ is {\it minimal} if it is minimal for the partial order $\prec$ on pairs of strings where $(v \sim v') \prec (w \sim w')$ if $v \prec w$ and $v' \prec w'$. 

Let $w_1$ and $w_2$ be any pair of shortlex words which satisfies $\overline{w_1 x} = \overline{w_2}$. Then $w_1 x$ and $w_2$ differ only by their minimal relation. It is shown that the size of such minimal relations are bounded (\cite[Lemma 4.3.2]{CE}), which completes the proof. \qed

\vspace*{\baselineskip}

From the proof of Theorem 2.2, the next lemma follows immediately. 

\vspace*{\baselineskip}

{\bf Lemma 2.3.} {\it Let $G$ and $\Sigma$ be the same as Theorem 2.2 and assume that $\Sigma$ is ordered. 
Let $w_1$ and $w_2$ be any pair of shortlex words which satisfies $\overline{w_1x} = \overline{w_2}$. Let $a_1^{r_1(t)} \cdots a_n^{r_n(t)}$ and $a_1^{s_1(t)} \cdots a_n^{s_n(t)}$ be the normal forms of $t$-prefixes of $w_1$ and $w_2$ respectively. Then there exists a constant $\kappa \ge 1$ such that $\max_t \{\sum_{i = 1}^n|r_i(t)-s_i(t)| \} \le \kappa$.}

\vspace*{\baselineskip}

In the followings, we regard $\kappa$ in Lemma 2.3 as the fellow travellers constant of a shortlex automatic structure for a finitely generated abelian group with an ordered generating set.

\section{Construction of word acceptors}

Let $G$ be a finitely generated abelian group and $\Sigma$ be an ordered symmetric generating set. There exists a shortlex automatic structure for $(G, \Sigma)$ from Theorem 2.2. In this section, we consider a canonical construction of word acceptors which are suited to compute the growth functions.

\begin{figure}
   \centerline{ 
    \epsfysize=8cm
   \epsfbox{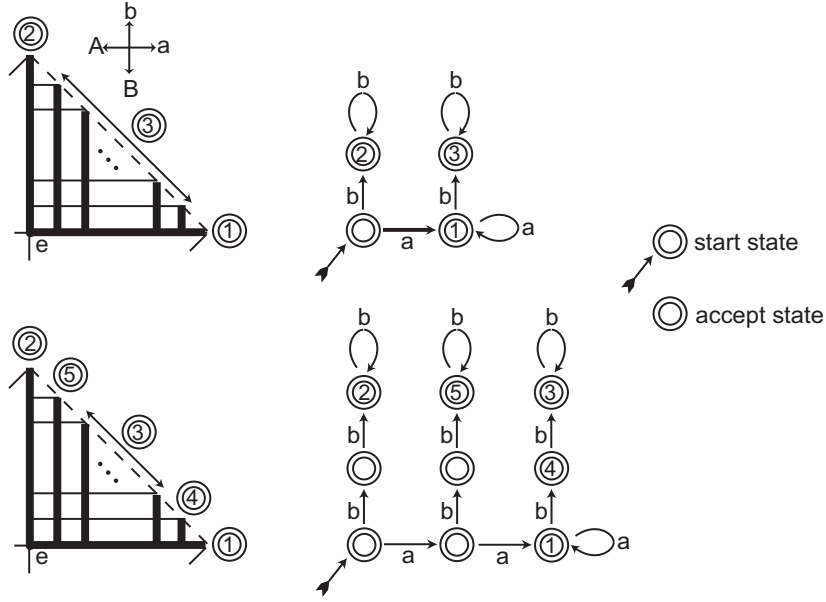}}
   \caption{Paths in Cayley graphs corresponding to word acceptors}
\end{figure}

First, we see a rough example which is stated in Figure 1. The left hand side in the figure depicts the Cayley graphs of $G = <a,b \ | \ ab =ba>$. The right hand side depicts word acceptors both of which accept the shortlex words for $(G, \Sigma)$ with the order $a <  b $. 
Arrows and states corresponding to the inverse elements and the failure states are omitted for simplicity. 
In the Cayley graphs, the dotted lines indicate $\partial \Gamma_n$. The thick lines indicate paths according to words accepted by the automata stated in the right. The numbers in circles in the figure indicate the correspondence between vertices of $\partial \Gamma_n$ and accept states of the word acceptors. 
Notice that the vertices of $\partial \Gamma_n$ are classified more minutely if we use the word acceptor below. 

Now we prepare some notations to describe the rigorous construction of the word acceptors. 
In an automaton, an {\it incoming }(resp. {\it outgoing}) {\it arrow} of a state $\sigma$ is an arrow whose target (resp. source) is $\sigma$. 
An arrow is called a {\it loop} if the source and target of the arrow coincide. In the followings, we use the notation ``incoming (resp. outgoing) arrows" only for arrows whose source and target states are distinct, thus loops are excluded from the set of incoming (resp. outgoing) arrows of any state.

We fix a presentation for $(G, \Sigma)$, index the relators and represent the $i$-th relator as $w_i = e$. Let $length(w_i)$ be the word length of $w_i$ and $\mu = \sum_i length(w_i)$. 
If a letter of $\Sigma$ may not occur infinitely in shortlex words for $(G, \Sigma)$, then the number of occurrences of the letter does not exceed $\mu$. 
Let $\gamma$ be any positive integer which is greater than $\mu$.

\vspace*{\baselineskip}

{\sc Definition} (canonical word acceptor). Construct a word acceptor which is part of a shortlex automatic structure for $(G, \Sigma)$ as follows.

\noindent
(1) Put the start state. 

\noindent
(2)  Put the unique failure state. Add $|\Sigma|$ loops, each of which is labelled by each letter of $\Sigma$, to the failure state.

\noindent
(3-1) For each state added in the previous procedure (1), (3-3a) or (3-3b), we continue the following steps. 

\noindent
(3-2) Let $\sigma$ be the added state and $\alpha$ be the label of the incoming arrow of $\sigma$. If $\sigma$ is the start state, we regard $\alpha$ as empty. 
Classify each letter of $\Sigma$ except $\alpha$ into three cases where the letter may occur infinitely, finitely or may not occur after the prefix words according to the transitions from the start state to $\sigma$.

\noindent
(3-3) Let $a_i$ be one of the classified letters in (3-2). In each case where the letter may occur infinitely, finitely or may not occur, do (3-3a), (3-3b) or (3-3c) respectively. 

\noindent
(3-3a) Construct a line which is rooted at $\sigma$ and which consists of arrows and states as the followings. 
First, we regard $\sigma$ as the end of the line whose length is equal to $0$. Add an arrow labelled by $a_i$ to the end, put its target state and regard the state as the new end. Continue these operations until the length of the line achieves $\gamma$. Finally add a loop labelled by $a_i$ to the end of the line. See Figure 2(a).

\noindent
(3-3b) Assume that $a_i$ may occur at most $\gamma_i$ times. Notice that $\gamma_i < \gamma$ by definition. Construct a line rooted at $\sigma$ whose length is equal to $\gamma_i$ the same as (3-3a). Finally, add an arrow, whose label is $a_i$ and whose target is the failure state, to the end of the line. See Figure 2(b).

\begin{figure}
   \centerline{ 
    \epsfysize=3cm
   \epsfbox{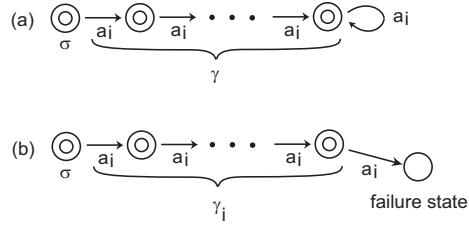}}
   \caption{Construction of word acceptor stated in (2-2)}
\end{figure}

\noindent
(3-3c) Add an arrow, whose label is $a_i$ and whose target is the failure state, to $\sigma$.

\noindent
(3-4) If there remains a letter classified in (3-2), we return to (3-3).

\noindent
(3-5) If there remains a state listed in (3-1), we return to (3-2). 

\noindent
(3-6) If there appear all letters of $\Sigma$ as labels of outgoing arrows or a loop in each state, we finish the construction. Otherwise, return to (3-1).

We regard all states except the failure state as accept states. 
The automaton constructed in this manner is called the $\gamma$-{\it canonical word acceptor} for $(G, \Sigma)$ and denoted by $W_{\gamma}$. 
\vspace*{\baselineskip}
\begin{figure}
   \centerline{ 
    \epsfysize=7cm
   \epsfbox{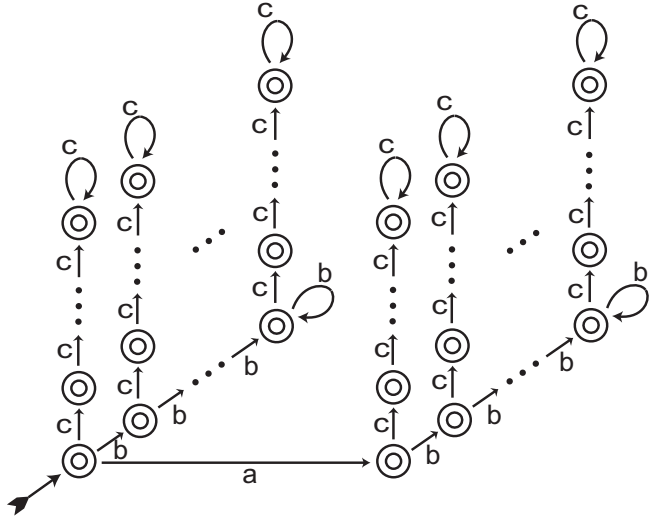}}
   \caption{The word acceptor $W_{\gamma}$ for the group stated in Example 3.1}
\end{figure}

Note that each state except the start state and the failure state has just one incoming arrow and at most one loop. Thus $W_{\gamma}$ is a tree if all loops, the failure state and all incoming arrows of the failure state are eliminated. Therefore, there exists exactly one path from the start state to each accept state in $W_{\gamma}$. 

\vspace*{\baselineskip}

{\sc Example 3.1.} Figure 3 depicts the $\gamma$-canonical automaton for some $\gamma$ in the case where $G = <a, b ,c \ | \ a^2 = b, \ ac = ca>$ with the order $a < b <  c $. Arrows and states corresponding to the inverse elements and the failure state are omitted. 
The lengths of oblique and vertical lines are equal to $\gamma$. 

\vspace*{\baselineskip}

In the followings, we often regard a finite state automaton as a labelled directed graph. A {\it walk} in a directed graph is a sequence of vertices $(v_0, v_1, \cdots, v_l)$ such that adjacent pairs $(v_0, v_1), (v_1, v_2), \cdots ,(v_{l-1}, v_l)$ are consistent with the directed edges of the graph including their directions. The number of edges in a walk is called the {\it length} of the walk. Notice that a path can be seen as a walk with no repeated vertices. In this paper, a walk whose length is equal to $l$ is called an $l$-{\it walk}.

Let $m$ be the number of accept states in an automaton and the accept states be numbered from $1$ to $m$. A row vector in $\{0,1\}^m$ is called the {\it characteristic vector} of the $k$-th state if the $k$-th entry is equal to $1$ and the others are equal to $0$. A $m \times m$ matrix is called the {\it transition matrix} of the automaton if each $(i,j)$-entry is the number of arrows from the $i$-th state to the $j$-th state. 

For the states of $W_{\gamma}$, we denote the start state by $\sigma_1$ and the other accept states by $\sigma_k$ for $2 \le k \le m$, where $m$ is the number of accept states of $W_{\gamma}$. 
Let $c_j(\sigma_k)$ be the number of $j$-walks from $\sigma_1$ to $\sigma_k$ in $W_{\gamma}$ and $C(\sigma_k, z) = \sum_{j = 0}^{\infty} c_j(\sigma_k) z^j$.

\vspace*{\baselineskip}

{\bf Lemma 3.2.} {\it For each state $\sigma_k$ of $W_{\gamma}$, $C(\sigma_k, z)$ is a rational function with integral coefficients. Furthermore the denominator of $C(\sigma_k, z)$ is $(1-z)^r$ with some nonnegative integer $r$ which does not exceed rank $G$ and the numerator is $z^q$ with some positive integer $q$. }

\vspace*{\baselineskip}

{\rm Proof.} Let $A$ be the transition matrix of $W_{\gamma}$ and ${\bf w}_i$ be the characteristic vector of $\sigma_i$ for $1 \le i \le m$. Then $C(\sigma_k,z) = \sum_{j=0}^{\infty}{\bf w}_1 (z A)^j {\bf w}_k^t$. Thus the proof consists in the computation of this formula.

First, we define a directed subgraph $W_{\gamma}^{(k)}$ of $W_{\gamma}$ as follows. The set of vertices of $W_{\gamma}^{(k)}$ consists of all states which are contained in the unique path from $\sigma_1$ to $\sigma_k$ in $W_{\gamma}$. The set of directed edges of $W_{\gamma}^{(k)}$ consists of all arrows which are contained in the path and all loops whose source and target states are contained in the set of vertices of $W_{\gamma}^{(k)}$. Then the number of $j$-walks from $\sigma_1$ to $\sigma_k$ in $W_{\gamma}$ coincides with that in $W_{\gamma}^{(k)}$.

Let $l$ be the length of the unique path from $\sigma_1$ to $\sigma_k$ in $W_{\gamma}^{(k)}$. In the case where $W^{(k)}_{\gamma}$ has no loops, $C(\sigma_k, z) = z^l$ and it satisfies the assertion of the statement. Thus it is sufficient to consider the case where $W^{(k)}_{\gamma}$ has loops. In that case, we refine $W^{(k)}_{\gamma}$ as follows. 
If some state $\sigma_i$ except $\sigma_1$ and $\sigma_k$ has no loop, we put an arrow from the previous state of $\sigma_i$ to the next state of that, and eliminate $\sigma_i$ with the incoming and outgoing arrows. If $\sigma_k$ has no loop, we eliminate $\sigma_k$ with the incoming arrow. Consequently, each of the remained states except $\sigma_1$ has a loop. 
Let $k'$ be the number of the remained states, thus the number of the eliminated ones is $l-k'+1$. We denote the remained graph by $V^{(k')}_{\gamma}$. Notice that there exists exactly one path in $V^{(k')}_{\gamma}$ which starts at $\sigma_1$ and whose length is equal to $k'-1$. Each state of $V^{(k')}_{\gamma}$ is renamed $\tau_i$, whose index is determined as the following. The state $\sigma_1$ of $W^{(k)}_{\gamma}$ turns into $\tau_1$. The other states are renumbered from $2$ to $k'$ according to the order of vertices along the unique path, thus the end of the path turns into $\tau_{k'}$. See Figure 4.

Let $B_{k'}$ be the transition matrix of $V^{(k')}_{\gamma}$. Then

\noindent

\[
B_{k'} = (b_{ij}) = \left\{
\begin{array}{l}
1 \ \ \ \ {\rm if} \ i = j \ \ \ {\rm for} \ 2 \le j \le k' \\
1 \ \ \ \ {\rm if} \ i = j-1 \ \ \ {\rm for} \ 2 \le j \le k' \\
0 \ \ \ \ {\rm otherwise.} 
\end{array}
\right.
\]
We denote by ${\bf v}_i$ the characteristic vector of $\tau_i$ in $V^{(k')}_{\gamma}$. Since the number of $j$-walks from $\sigma_1$ to $\sigma_k$ in $W^{(k)}_{\gamma}$ coincides with the number of $\{j-(l-k'+1)\}$-walks from $\tau_1$ to $\tau_{k'}$ in $V^{(k')}_{\gamma}$, we have

\[
C(\sigma_k,z) = \sum_{j=0}^{\infty}{\bf w}_1 (z A)^j {\bf w}_k^t = \sum_{j=0}^{\infty}{\bf v}_1 z^{l-k'+1} (z B_{k'})^j {\bf v}_{k'}^t \]
\[
 ={\bf v}_1 z^{l-k'+1} (I -zB_{k'})^{-1} {\bf v}_{k'}^t=\frac{z^{l}}{(1-z)^{k'-1}}.
\]

\noindent
It is obvious from the fundamental theorem of abelian groups that $k'-1$ does not exceed rank $G$. \qed

\vspace*{\baselineskip}

For the proof of the main theorem, we prepare a general formula derived from Lemma 3.2. In the next lemma, we follow the notations in the proof of Lemma 3.2.

\vspace*{\baselineskip}

{\bf Lemma 3.3.} {\it Let $\eta_1$ and $\eta_2$ be positive integers. Then $\sum_{j=\eta_1}^{\infty}{\bf w}_1 z^{\eta_2} (z A)^j {\bf w}_k^t$ is a sum of rational functions each of whose denominators is $(1-z)^{i}$ for $0 \le i \le k'-1$ and each of whose numerators is $p_i z^{q_i}$ such that each $q_i$ is a positive integer, each $p_i$ is a nonnegative integer for $0 \le i \le k'-2$ and $p_{k'-1} =1$.}

\vspace*{\baselineskip}

{\rm Proof.} Since 

\[ \sum_{j=\eta_1}^{\infty}{\bf w}_1 z^{\eta_2} (z A)^j {\bf w}_k^t =\sum_{j=\eta_1}^{\infty}{\bf v}_1 z^{\eta_2+l-k'+1} (z B_{k'})^j {\bf v}_{k'}^t \]

\[=  \sum_{j=0}^{\infty}{\bf v}_1 z^{\eta_1+\eta_2+(l-k'+1)} B_{k'}^{\eta_1} (z B_{k'})^j {\bf v}_{k'}^t
\]

\noindent
and the $(1,2)$-entry of $B_{k'}^{\eta_1}$ for any $\eta_1 \ge 1$ is equal to $1$, the statement follows immediately. \qed

\begin{figure}
   \centerline{ 
    \epsfysize=4.5cm
   \epsfbox{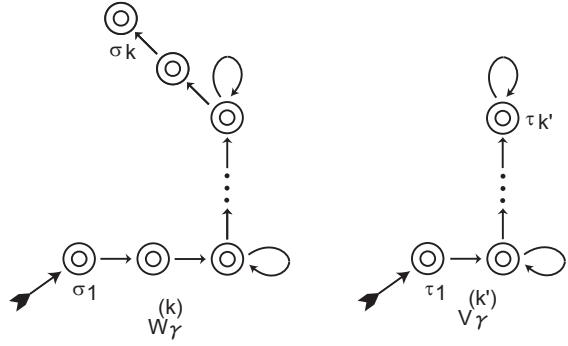}}
   \caption{Construction of $W^{(k)}_{\gamma}$ and $V^{(k')}_{\gamma}$}
\end{figure}

\section{Main argument}

In this section, we prove the main theorem. 
Let $G$, $\Sigma$ and $S$ be the same as the statement of the main theorem. Let $d$ be the diameter of $S$. Fix a base point $p \in S$. 
Let $\kappa$ be the fellow travellers constant of a shortlex automatic structure for $(G, \Sigma)$ defined in the paragraph after Lemma 2.3. We fix a positive integer $\gamma$ which is greater than $d\kappa + \mu$ and construct $W_{\gamma}$ for $(G, \Sigma)$. 

Let $v$ be a vertex of $\Gamma(G,\Sigma)$. To simplify notation, we use the same symbol $v$ for the corresponding element of $G$. A word $w$ is called an {\it accepted word} of $v$ if $w$ is accepted by $W_{\gamma}$ and $\overline{w} = v$. Notice that there exists exactly one accepted word for an element of $G$ from the uniqueness property of shortlex words. A path from the identity element to $v$ in $\Gamma$ is called the {\it accepted path} of $v$ if the sequence of the labels along the path coincides with the accepted word of $v$. 
The accepted paths are geodesics in $\Gamma$ since $W_{\gamma}$ accepts only shortest words among representatives for elements of $G$. 
When we follow the arrows in $W_{\gamma}$ according to the labels of the accepted word of $v$ starting at $\sigma_1$, we end at some state of $W_{\gamma}$. The state is called the {\it final state} of $v$.
In the followings, it is necessary to pay attention to the difference between walks, including paths, in $W_{\gamma}$ and accepted paths in $\Gamma$.

\vspace*{\baselineskip}
{\bf Lemma 4.1.} {\it Assume that $n$ is greater than $\gamma | \Sigma |+d$. Let $f : S \to \Gamma_n$ be a morphism such that $fS \not\subset \Gamma_{n-1}$ and $\sigma_k$ be the final state of $fp$. Let $f'$ be a morphism from $S$ to $\Gamma$ such that the final state of $f'p$ is $\sigma_k$. Then $d_{\Gamma}(e,f'p) = d_{\Gamma}(e, fp)$ if and only if $f'S \subset \Gamma_n$ and $f'S \not\subset \Gamma_{n-1}$.}

\vspace*{\baselineskip}

{\rm Proof.} Let $a_1^{r_1} \cdots a_l^{r_l}$ be the normal form of the accepted word of $fp$. 
From the assumption that $n$ is greater than $\gamma | \Sigma |+d$, there exist indices $i_1, \cdots, i_m$ such that $r_{i_1}, \cdots, r_{i_m}$ are greater than or equal to $\gamma$. 
Thus the unique path from $\sigma_1$ to $\sigma_k$ in $W_{\gamma}$ contains exactly $m$ states each of which is the source and target of a loop. See Figure 5. 
Let $a_1^{r'_1} \cdots a_l^{r'_l}$ be the normal form of the accepted word of $f'p$. 
Since the final state of $f'p$ is $\sigma_k$, the walk according to $a_1^{r'_1} \cdots a_l^{r'_l}$ in $W_{\gamma}$ contains the path from $\sigma_1$ to $\sigma_k$. Thus $r'_i = r_i$ for any $i \neq i_1, \cdots, i_m$ and $r'_{i_j} \ge \gamma$ for $j = 1, \cdots, m$. 

\begin{figure}
   \centerline{ 
    \epsfysize=3cm
   \epsfbox{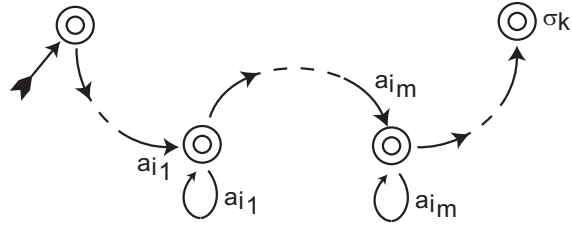}}
   \caption{The unique path in $W_{\gamma}$}
\end{figure}

Let $q$ be any vertex of $S$ and $a_1^{s_1} \cdots a_l^{s_l}$ be the normal form of the accepted word of $fq$. By Lemma 2.3, $|r_{i_j}-s_{i_j}| \le d\kappa$, thus $s_{i_j} \ge \mu$ for $j = 1, \cdots, m$. Consider a word $a_1^{s'_1} \cdots a_l^{s'_l}$ such that $s'_i = s_i$ for any $i \neq i_1, \cdots i_m$ and $s'_{i_j} = s_{i_j}+ (r'_{i_j} - r_{i_j})$ for $j = 1, \cdots, m$. Since $r'_{i_j} \ge \gamma$ and $|r_{i_j}-s_{i_j} | \le d\kappa$, $s'_{i_j} \ge \mu$ for $j = 1, \cdots, m$. Combining this observation with the assumption that $a_1^{s_1} \cdots a_l^{s_l}$ is shortlex, we conclude that $a_1^{s'_1} \cdots a_l^{s'_l}$ is also the normal form of a shortlex word. Furthermore, $a_1^{s'_1} \cdots a_l^{s'_l}$ is a representative for $f'q$ since $\Gamma$ is homogeneous and $G$ is commutative. Therefore $a_1^{s'_1} \cdots a_l^{s'_l}$ is the normal form of the accepted word of $f'q$.

Assume that $d_{\Gamma}(e,f'p) = d_{\Gamma}(e, fp)+ \eta$, it means that $\sum_{i = 1}^{l} r'_i = \sum_{i = 1}^{l} r_i + \eta$, for some integer $\eta$. Then $\sum_{i =1}^l s'_{i} = \sum_{i=1}^l s_{i}+ \eta$ from the construction of $s'_i$, and it follows that $d_{\Gamma}(e,f'q) = d_{\Gamma}(e, fq)+ \eta$. Thus the statement follows immediately. \qed

\vspace*{\baselineskip}

Let $n$, $f$ and $\sigma_k$ be the same as the assumptions of Lemma 4.1. Let $d_{\Gamma}(e,fp) = n-\delta_k$, thus $\delta_k$ satisfies $0 \le \delta_k \le d$. 
Each $(n-\delta_k)$-walk from $\sigma_1$ to $\sigma_k$ in $W_{\gamma}$ determines a unique morphism $f': S \to \Gamma$ because the natural map $\pi : L(W_{\gamma}) \to G$ is bijective. Furthermore $f'$ satisfies $f'S \subset \Gamma_n$ and $f'S \not\subset \Gamma_{n-1}$ from Lemma 4.1. On the contrary, each morphism $f':S \to \Gamma_n$ such that $f'S \not\subset \Gamma_{n-1}$ and that the final state of $f'p$ is $\sigma_k$ determines a unique $(n-\delta_k)$-walk from $\sigma_1$ to $\sigma_k$ in $W_{\gamma}$ according to the accepted word of $f'p$ because of the uniqueness of shortlex words. 
Thus there is a one to one correspondence between the set of walks from $\sigma_1$ to $\sigma_k$ in $W_{\gamma}$ whose lengths are $n-\delta_k$ and the set of morphisms from $S$ to $\Gamma_n$ such that the images of $S$ are not contained in $\Gamma_{n-1}$ and that the final states of the images of $p$ are $\sigma_k$.

\vspace*{\baselineskip}

{\rm Proof of the main theorem.} Let $c_n(\sigma_k, S)$ be the number of morphisms from $S$ to $\Gamma_n$ such that the images of $S$ are not contained in $\Gamma_{n-1}$ and that the final states of the images of $p$ coincide with $\sigma_k$. Let $C(\sigma_k, S, z) = \Sigma_{j=0}^{\infty} c_j(\sigma_k, S) z^j$, then $C(S,z) = \Sigma_k C(\sigma_k, S, z)$. Let $c_j(\sigma_k)$ be the number of walks from $\sigma_1$ to $\sigma_k$ in $W_{\gamma}$ whose lengths are equal to $j$ as stated in section 3.

From the paragraph after Lemma 4.1, $c_j(\sigma_k, S)$ is equal to $c_{j-\delta_k}(\sigma_k)$ for $j > \gamma|\Sigma| + d$. Let $A$ be the transition matrix of $W_{\gamma}$ and ${\bf w}_i$ be the characteristic vector of $\sigma_i$. 
Then we have

\[
C(\sigma_k, S, z) =\sum_{j =0}^{\infty} c_j(\sigma_k, S) z^j= \sum_{j =0}^{\gamma | \Sigma| + d} c_j(\sigma_k, S) z^j+\sum_{j =\gamma|\Sigma|+d+1}^{\infty} c_{j-\delta_k}(\sigma_k) z^{j}
\]
\[
 = \sum_{j =0}^{\gamma |\Sigma|+ d} c_j(\sigma_k,S) z^j + \sum_{j =\gamma |\Sigma|+d+1}^{\infty} {\bf w}_1 (zA)^{j-\delta_k} {\bf w}_k^t \]

\[= \sum_{j =0}^{\gamma |\Sigma|+d} c_j(\sigma_k,S) z^j + \sum_{i=0}^{r_k} \frac{p^{(k)}_{i} z^{q^{(k)}_{i}}}{(1-z)^{i}}
\]

\noindent
such that $r_k$ is a nonnegative integer, each $q^{(k)}_{i}$ is a positive integer, each $p^{(k)}_{i}$ for $0 \le i \le r_k-1$ is a nonnegative integer and $p^{(k)}_{r_k} = 1$ from Lemma 3.3. Further, each $r_k$ does not exceed $rank \ G$ for any $k$ and there exists an index $k$ such that $r_k = rank \ G$ from the fundamental theorem for abelian groups. 

Let $r = rank \ G$. Then we have

\[
C(S, z) = \sum_k C(\sigma_k, S, z) = \frac{(1-z)P_0(z)+ P_1(z)}{(1-z)^r}
\]

\noindent
such that $P_0(z)$ and $P_1(z)$ are polynomials with integral coefficients and that all coefficients of $P_1(z)$ are positive. 
Since $P_1(1) > 0$, $(1-z)$ is not a factor of the numerator of the last rational function. Consequently, the last rational function is irreducible. \qed

\section{Examples}

In this section, we see two examples. The first one is presented in~\cite{E}. Recall that $b_n(S)$ denotes the number of morphisms $f:S \to \Gamma_n$.

\vspace*{\baselineskip}

{\sc Example 5.1} (\cite{E}). Let $G = <a , b, c \ | \ ab = ba, ac = ca, bc = cb>$. The growth functions for the subgraphs $S = \{ vertex \}, \{ a \}, \{ b \}, \{ c \}, \{a,b\}$ and $\{a,b,c\}$ are computed as follows.

\[
C(\{vertex\}, z) = 1+7z+18z^2+38z^3 + \cdots = \frac{(z+1)^3}{(1-z)^3}.
\]
\[
C(\{ a \}, z) = C(\{ b \}, z) =C(\{ c \}, z) = 2z + 10 z^2 + 26z^3+50z^4 + \cdots = \frac{2z(z+1)^2}{(1-z)^3}.
\]
\[
C(\{a,b\}, z) = z+7z^2+21z^3 +43z^4+ \cdots = \frac{z(z+1)(3z+1)}{(1-z)^3}.
\]
\[
C(\{a,b,c\}, z) = 4z^2+16z^3 + 36z^4+ 64z^5+\cdots = \frac{4z^2(z+1)}{(1-z)^3}.
\]

In~\cite{E}, $B(S,z) = \sum_{j=0}^{\infty} b_j(S) z^j$ are stated instead of $C(S,z)$. It is clear that the functions $C(S,z)$ are computed by the formula $C(S,z) = B(S,z) - zB(S,z)$. 

\vspace*{\baselineskip}

{\sc Example 5.2.} Let $G =<a , b, c \ | \ ab = ba, c = ab>$ thus the rank of $G$ is equal to $2$. The growth functions for the subgraphs $S = \{ vertex \}, \{ a \}, \{ b \}, \{ c \}$ and $\{a,b,c\}$ are computed as follows.

\[
C(\{vertex\}, z) = 1+ 6z + 12z^2 + 18z^3 + \cdots = \frac{1+4z+z^2}{(1-z)^2}.
\]
\[
C(\{ a \}, z) = C(\{ b \}, z) =C(\{ c \}, z) =2z + 8z^2 + 14z^3 + 20z^4 + \cdots = \frac{2z(1+2z)}{(1-z)^2}.
\]
\[
C(\{a,b,c\}, z) =  z + 6z^2+ 12z^3 + 18z^4 + \cdots = \frac{z+4z^2+z^3}{(1-z)^2.}
\]

\vspace*{\baselineskip}

\end{document}